\documentclass[10pt,twoside,a4paper,reqno]{amsart}
\usepackage{amscd,amsmath,amsthm,amsfonts,latexsym,amssymb}

\theoremstyle{plain}
\newtheorem{theo+}           {Theorem}
\newtheorem{prop+}           {Proposition}
\newtheorem{coro+}           {Corollary}
\newtheorem{lemm+}           {Lemma}
\newtheorem{conjecture}      {Conjecture}

\theoremstyle{definition}
\newtheorem{rema+}           {Remark}
\newtheorem{defi+}           {Definition}
\newtheorem{problem}         {Problem}
\newtheorem*{ack}            {Acknowledgement}
\newtheorem{nota}            {Notation}

\newenvironment{theorem}{\begin{theo+}}{\end{theo+}}
\newenvironment{proposition}{\begin{prop+}}{\end{prop+}}
\newenvironment{corollary}{\begin{coro+}}{\end{coro+}}
\newenvironment{lemma}{\begin{lemm+}}{\end{lemm+}}
\newenvironment{remark}{\begin{rema+}}{\end{rema+}}
\newenvironment{definition}{\begin{defi+}}{\end{defi+}}

\newcommand {\bC} {\mathbb {C}}
\newcommand {\bR} {\mathbb {R}}
\newcommand {\bN} {\mathbb {N}}
\newcommand {\bZ} {\mathbb {Z}}
\newcommand {\bs} {\mathbf {s}}
\newcommand {\bt} {\mathbf {t}}

\newcommand {\al} {\alpha}

\newcommand {\de} {\delta}
\newcommand {\la} {\lambda}
\newcommand {\ga} {\gamma}
\newcommand {\La} {\Lambda}
\newcommand {\si} {\sigma}

\newcommand {\Ga} {\Gamma}
\newcommand {\vf} {\varphi}
\newcommand {\ve} {\varepsilon}
\newcommand {\pc} {\preccurlyeq}
\newcommand {\pr} {\prec}
\newcommand {\calLP} {\mathcal {LP}}
\newcommand {\calPS} {\mathcal {PS}}
\newcommand {\calH} {\mathcal {H}}

\newcommand {\calA} {\mathcal {A}}

\newcommand {\calZ} {\mathcal {Z}}
\newcommand {\calD} {\mathcal {D}}
\newcommand {\calT} {\mathcal {T}}
\newcommand {\calX} {\mathcal {X}}
\newcommand {\calF} {\mathcal {F}}
\newcommand {\vP} {\varPi}
\newcommand {\E} {\text{End}\,}
\newcommand {\li} {l^{\infty}}
\newcommand {\pa} {\partial}
\newcommand {\fm} {\mathfrak{m}}

\begin{document}

\numberwithin{equation}{section}

\title[Spectral order and isotonic differential operators]{Spectral order and 
isotonic differential operators of Laguerre-P\'olya type}
\author[J.~Borcea]{Julius Borcea}
\address{Department of Mathematics, Stockholm University, SE-106 91 Stockholm,
   Sweden}
\email{julius@math.su.se}



\begin{abstract}
The spectral order on $\bR^n$ induces a natural partial 
ordering on the manifold $\calH_{n}$ of monic hyperbolic polynomials of 
degree $n$. We show that all differential operators of Laguerre-P\'olya type 
preserve the spectral order. We also establish a global monotony property for 
infinite families of deformations of these operators parametrized by the 
space $\li$ of real bounded sequences. As a consequence, we deduce that the 
monoid $\calA'$ of linear operators that preserve averages of zero sets and 
hyperbolicity consists only of differential operators of 
Laguerre-P\'olya type which are both extensive and isotonic. In particular, 
these results imply that any hyperbolic polynomial is the global minimum of 
its $\calA'$-orbit and that Appell polynomials are characterized by a 
global minimum property with respect to the spectral order.
\end{abstract}

\maketitle

\section*{Introduction and main results}\label{s0}

This is the third part of a series of papers \cite{B,BBS,BP,BS} on the 
connections between linear operators acting on partially ordered manifolds of 
polynomials, the distribution of zeros of polynomials, and the theory of 
majorization. 

Linear differential operators acting on various function spaces and 
classical majorization have both been extensively studied albeit so far only 
in separate contexts. On the one hand, differential operators of 
infinite order appear naturally in many applications. From a topological 
point of view they form a total set of linear continuous operators between 
spaces of differentiable functions \cite{K}, which is rather reminiscent of 
Peetre's abstract characterization of differential operators \cite{P}. In this
paper we are mainly concerned with linear operators of 
Laguerre-P\'olya type, that is, infinite order differential operators 
induced by the Laguerre-P\'olya class 
of entire functions. The significance of the latter 
stems from the fact that 
it consists precisely of those functions 
which are locally uniform limits in $\bC$ of sequences of polynomials with 
all real zeros \cite{L}. There is a very rich literature on differential 
operators of 
Laguerre-P\'olya type and their applications to the study of the distribution 
of zeros of certain Fourier transforms, P\'olya-Schoenberg frequency functions 
and totally positive matrices, the inversion and representation theories of 
convolution transforms, and the final set problem for trigonometric 
polynomials. Recently, such operators were also studied in connection with 
various generalizations of the P\'olya-Wiman conjecture. Further details 
on these topics and related questions may be found in 
e.g.~\cite{CC1,CC2,KOW} and references therein.

On the other hand, the notion of (classical) majorization was first studied
by economists early in the twentieth century as a means for altering the 
unevenness of distribution of wealth or income. Classical majorization was a 
key tool in Schur's work on Hadamard's determinantal inequality and 
the spectra of positive semide\-finite Hermitian matrices \cite{DK}. This 
notion was later formalized as a preorder on $n$-vectors of 
real numbers -- also known as the spectral order on $\bR^n$ -- by Hardy, 
Littlewood and P\'olya in their study of symmetric means and analytic 
inequalities \cite{HLP}. The spectral order has since found important 
applications in operator theory, convex analysis,
combinatorics and statistics \cite{An1,An2,MO}. 
As recent results have shown, classical majorization plays also a remarkable 
role in the study of quantum state mixing and efficient 
measurements in quantum mechanics \cite{NV}, quantum algorithm 
design \cite{LMD} and the analysis of entanglement transformations in quantum 
computation and information theory \cite{JP}.

As we explain below, the spectral order on $\bR^n$ induces a natural partial 
ordering $\pc$ on the 
manifold $\calH_n$ of monic univariate polynomials of degree $n$ with all real 
zeros (cf.~\cite{B,BS}). Polynomials of this type are often called 
{\em hyperbolic} owing to the standard terminology used in the theory of 
partial differential equations \cite{G}, singularity theory and related topics
\cite{Ar}. Let $\vP:=\bC[x]$ be the space of 
complex univariate polynomials 
regarded as functions on the complex plane. The main purpose of this paper is 
to study the properties of the 
posets $(\calH_n,\pc)$, $n\in \bN$, under the action of 
hyperbolicity-preserving linear operators, that is, operators acting on $\vP$ 
that map hyperbolic polynomials to 
hyperbolic polynomials. 
Given a monic polynomial $P\in \vP$ with 
$\deg P=n\ge 1$ we define $\calZ(P)$ to be the unordered $n$-tuple consisting 
of the zeros of $P$, each zero occurring as many times as its multiplicity. 
Thus $\calZ(P)\in \bC^n/\Sigma_n$, where $\Sigma_n$ is the symmetric group on 
$n$ elements. We denote by $\Re \calZ(P)$ the unordered $n$-tuple whose 
components are the real parts of the points in $\calZ(P)$. Note that $P$ 
is hyperbolic if and only if $\Re \calZ(P)=\calZ(P)$. A hyperbolic polynomial 
with simple zeros is called 
{\em strictly hyperbolic}. Let $\calH_{n}\subset \vP$ be the real manifold of 
monic hyperbolic polynomials of degree $n$. We extend this notation to $n=0$ 
by setting $\calH_0=\{1\}\subset \vP$. Clearly, for $n\ge 1$ one has a natural 
set-theoretic identification between $\calH_n$ and $\bR^n/\Sigma_n$ by means 
of the {\em root map} 
\begin{equation}\label{rmap}
\begin{split}
\calZ:\calH_n & \longrightarrow \bR^n/\Sigma_n\\
P & \longmapsto \calZ(P).
\end{split}
\end{equation}

The following theorem is due to Hardy, Littlewood and P\'olya~\cite{HLP}:

\begin{theorem}\label{spec}
Let $X=(x_{1},x_{2},\ldots,x_{n})\in \bR^n/\Sigma_n$, 
$Y=(y_{1},y_{2},\ldots,y_{n})\in \bR^n/\Sigma_n$. The following conditions 
are equivalent:
\begin{enumerate}
\item[(i)] For any convex function $f: \bR \to \bR$ one
has $\sum_{i=1}^{n}f(x_{i})\le \sum_{i=1}^{n}f(y_{i})$.
\item[(ii)] There exists a doubly stochastic $n\times n$ matrix $A$ such that
$\tilde{X}=A\tilde{Y}$, where $\tilde{X}$ and $\tilde{Y}$ are column 
$n$-vectors obtained by some (and then any) ordering of the 
components of $X$ and $Y$, respectively.
\end{enumerate}
\end{theorem}

Theorem~\ref{spec} defines what is usually known as classical 
majorization or the spectral order on $\bR^n$: if the conditions of the 
theorem are satisfied we say that $X$ is {\em majorized} by $Y$ or that $X$ 
is {\em less than} $Y$ {\em in the spectral order}, which we denote by 
$X\pr Y$. One can easily check that if $X\pr Y$ then
$\sum_{i=1}^nx_{i}=\sum_{i=1}^ny_{i}$. Note that although 
the spectral order is only a preordering on $\bR^n$, Birkhoff's theorem 
\cite[Theorem 2.A.2]{MO} implies that it actually induces a partial ordering 
on $\bR^n/\Sigma_n$. Therefore, Theorem~\ref{spec} allows us to define a 
poset structure $(\calH_n,\pc)$ by setting $Q\pc P$ 
whenever $P,Q\in \calH_{n}$ and $\calZ(Q)\pr \calZ(P)$. In this way we may 
view the spectral order on $\bR^n$ as a natural partial ordering on the 
manifold $\calH_{n}$, which we call the {\em spectral order on} $\calH_{n}$.

We can now state the following isotonicity theorem, which is our first main 
result:

\begin{theorem}\label{main1} 
Let $n\ge 1$ and $P,Q\in \calH_n$ be such that $Q\pc P$. Then for any 
$\la\in \bR$ one has $Q-\la Q'\pc P-\la P'$. 
\end{theorem}

This has several natural consequences. Recall that by a 
classical result of P\'olya all 
differential operators of Laguerre-P\'olya type are hyperbolicity-preserving, 
see e.g.~\cite[Theorem 5.4.13]{RS}. Theorem~\ref{main1} implies that 
much more is actually true, namely all such operators 
preserve in fact the spectral order (Corollary~\ref{iso}). In 
particular, any degree-preserving differential operator of Laguerre-P\'olya 
type is isotonic with respect to the partial ordering $\pc$ on the 
manifold $\calH_n$ for all $n\in \bN$ (Corollary~\ref{coroiso1}). This gives 
a new characterization of the sequence of Appell polynomials
associated with an arbitrary function in the Laguerre-P\'olya class by means 
of a global minimum property with respect to the spectral order 
(Corollary~\ref{coroiso2}). 

Let $D=d/dx$ denote differentiation with respect to $x$. The second main 
result of this paper is the following monotonicity theorem:

\begin{theorem}\label{main2}
Fix $n\ge 1$ and let $\la_1,\la_2\in\bR$ be such that $\la_1 \la_2\ge 0$ and 
$|\la_1|\le |\la_2|$. Then $(1-\la_1D)e^{\la_1D}P\pc (1-\la_2D)e^{\la_2D}P$
for any $P\in\calH_n$.
\end{theorem}

Theorem~\ref{main2} allows us to study the orbit of an arbitrarily given 
hyperbolic polynomial
under the action of the monoid of differential operators of Laguerre-P\'olya 
type. We equip the space $\li$ of real bounded sequences with a natural 
partial ordering $\leqslant$ and define infinite families of deformations of 
differential operators of Laguerre-P\'olya type which are parametrized by 
vectors in $\li$. From Theorem~\ref{main2} we deduce that any such 
family satisfies a global monotony property with respect to both 
partial orderings $\pc$ and $\leqslant$ (Corollary~\ref{monot}). Moreover, 
these partial orderings are compatible with each other 
(Corollary~\ref{coromonot3}). It follows that the monoid 
$\calA'$ of all linear operators that act on each of the manifolds $\calH_n$,
$n\ge 1$, and 
preserve averages of zero sets consists only of differential operators of 
Laguerre-P\'olya type which are extensive with respect to $\pc$
(Corollary~\ref{coromonot1}). Thus, any hyperbolic polynomial is the global 
minimum of its $\calA'$-orbit with respect to the spectral order 
(Corollary~\ref{coromonot2}). 

The above results have further 
applications to the distribution of zeros of hyperbolic polynomials under 
the action of differential operators of Laguerre-P\'olya type 
(Corollaries~\ref{con1}-\ref{con3}). At the same time, they seem to suggest 
even deeper connections between linear (differential) operators, the 
distribution of zeros of real entire functions, and the theory of majorization.
As we point out in \S \ref{s4}, it would be interesting to 
know whether appropriate modifications of the aforementioned results could 
hold for transcendental entire functions in the Laguerre-P\'olya class. On the 
other hand, these results and those of \cite{B,BP,BS} hint at the possible 
existence of an ``analytic theory of classical majorization'' and may 
therefore also be seen as natural steps towards developing such a theory. 
Problem~2 in \cite{B} and Problems~\ref{pb1}-\ref{pb3} in \S \ref{s4} 
are intended as further steps in this direction.

\begin{ack}
The author would like to thank the anonymous referee for many useful 
suggestions and remarks. 
\end{ack}

\section{Theorem~\ref{main1} and applications}\label{s2}

\subsection{Proof of Theorem~\ref{main1}}\label{s21}

A key ingredient in the proofs of Theorems~\ref{main1} and~\ref{main2} is the 
following criterion due to Hardy, Littlewood and P\'olya \cite{HLP}.

\begin{theorem}\label{crit}
Let $X=(x_{1}\le x_{2} \le\ldots \le x_{n})$ and
$Y=(y_{1}\le y_{2} \le \ldots \le y_{n})$ be two $n$-tuples of 
real numbers. Then $X\prec Y$ if and only if the $x_{i}$'s and the 
$y_{i}$'s 
satisfy the following conditions:
$$\sum_{i=1}^{n}x_{i}=\sum_{i=1}^{n}y_{i}\,\text{ and }\,
\sum_{i=0}^{k}x_{n-i}\le \sum_{i=0}^{k}y_{n-i}\,\text{ for }\,
0\le k\le n-2.$$
\end{theorem}

We also make extensive use of {\em contractions}, a special kind of 
degree-preserving transformations acting on hyperbolic polynomials 
that we define as follows.

\begin{definition}\label{defcontr}
Let $P(x)=\prod_{i=1}^{n}(x-x_{i})\in \calH_{n}$, $n\ge 2$, and 
$1\le k<l\le n$. 
Assume that $x_{i}\le x_{i+1}$, $1\le i\le n-1$, and that $x_{k}\neq x_{l}$. 
Let further $t\in \left(0,\frac{x_{l}-x_{k}}{2}\right]$ and define 
$Q\in \calH_{n}$ to be 
the polynomial with zeros $y_{i}$, $1\le i\le n$, where $y_{k}=x_{k}+t$, 
$y_{l}=x_{l}-t$, and $y_{i}=x_{i}$, $i\neq k,l$. The polynomial $Q$ is 
called the {\em contraction of} $P$ {\em of type} $(k,l)$ 
{\em and coefficient} $t$ and is denoted by $Q=\calT(k,l;t)P$. The 
contraction $\calT(k,l;t)$ of $P$ is called {\em simple} if $l=k+1$ 
and it is called {\em nondegenerate} if $t\neq \frac{x_{l}-x_{k}}{2}$.
\end{definition}

\begin{remark}\label{T-tr}
The simple nondegenerate contractions in Definition~\ref{defcontr} may be 
viewed as elementary versions of the so-called $T$-transforms for $n$-tuples 
of real numbers. The latter are essentially a mathematical formulation of 
Dalton's ``principle of transfers'' (see~\cite{MO}) and 
were first used by Hardy, Littlewood and P\'olya in~\cite{HLP}.
\end{remark}

The proof of Theorem~\ref{main1} builds on several auxiliary technical 
results. The first two of these, Proposition~\ref{contr1} and Lemma~\ref{c1}
below, may also be restated in terms of $n$-tuples of real numbers or doubly 
stochastic matrices in view of Theorems~\ref{spec} and~\ref{crit}. However, 
since we are interested in the dynamics of polynomial zeros under the action 
of certain operators, it is convenient to formulate all the results 
exclusively in terms of polynomials.

\begin{proposition}\label{contr1}
Let $P,Q\in \calH_{n}$ be two distinct strictly hyperbolic polynomials such 
that $Q\pc P$. Then there exists a finite sequence of strictly hyperbolic 
polynomials $P_{1},\ldots,P_{m}\in \calH_{n}$ such that $P_{1}=P$, $P_{m}=Q$ 
and $P_{i+1}$ is a simple nondegenerate contraction of $P_{i}$ for 
$1\le i\le m-1$.
\end{proposition}

The algorithm described in the next lemma will be used to give a constructive 
proof of Proposition~\ref{contr1}. 

\begin{lemma}\label{c1}
Let $a<b$, $\si\in \left(0,\frac{b-a}{2}\right)$ and $p\in \bN$. Assume that 
$z_i$, $1\le i\le p$, are real numbers that satisfy 
$a+\si<z_1<\ldots<z_p<b-\si$ and set
$$P(x)=(x-a)(x-b)\prod_{i=1}^{p}(x-z_i),\,\,
Q(x)=(x-a-\si)(x-b+\si)\prod_{i=1}^{p}(x-z_i).$$ 
There exist simple nondegenerate contractions $\calT_1,\ldots,\calT_s$ such 
that $Q=\calT_s\cdots\calT_1 P$.
\end{lemma}

\begin{proof}
Set $x_1=a$, $x_{p+2}=b$ and $x_i=z_{i-1}$,
$2\le i\le p+1$, so that we may write 
$$P(x)=\prod_{i=1}^{p+2}(x-x_i) \text{ with } x_i<x_{i+1},\,1\le i\le p+1.$$
Choose $d\in \bN$ such that $\si<2^{d-1}\min(z_1-a-\si,\,b-z_p-\si)$ if $p=1$ 
and
$$\si<2^{d-1}\min\left(z_1-a-\si,\,b-z_p-\si,
\min_{1\le i\le p-1}(z_{i+1}-z_i)\right)\text{ if } p\ge 2.$$
We let $t=\frac{\si}{2^d}$ and build a finite sequence of polynomials 
$\{S_{1,i}\}_{i=0}^{p+1}$ as follows: 
$$S_{1,0}=P\text{ and }S_{1,i}=\calT(i,i+1;t)S_{1,i-1},\, 1\le i\le p+1.$$
Clearly, the contractions used in constructing this sequence 
are all simple. These contractions are also nondegenerate since
$$x_{i+1}-(x_i-t)>2t,\,1\le i\le p+1,$$
by the choice of $t$. Thus, all polynomials $S_{1,i}$, $0\le i\le p+1$, are 
strictly hyperbolic. In particular, this is true for the polynomial
$$P_1(x):=S_{1,p+1}(x)=\prod_{i=1}^{p+2}\left(x-x_{i}^{(1)}\right),$$
where $x_{1}^{(1)}=a+t$, $x_{p+2}^{(1)}=b-t$ and $x_{i}^{(1)}=z_{i-1}$,
$2\le i\le p+1$, so that $x_{i}^{(1)}<x_{i+1}^{(1)}$ for $1\le i\le p+1$. We 
now use the same contractions as above to construct a finite sequence of 
polynomials $\{S_{2,i}\}_{i=0}^{p+1}$ starting with the polynomial $P_1$:
$$S_{2,0}=P_1\text{ and }S_{2,i}=\calT(i,i+1;t)S_{2,i-1},\, 1\le i\le p+1.$$
Repeating this procedure $r$ times we arrive at the polynomial
$$P_r(x):=S_{r,p+1}(x)=\prod_{i=1}^{p+2}\left(x-x_{i}^{(r)}\right),$$
where $x_{1}^{(r)}=a+rt$, $x_{p+2}^{(r)}=b-rt$ and $x_{i}^{(r)}=z_{i-1}$ for
$2\le i\le p+1$. It is clear that all the contractions used in constructing
the polynomial $P_r$ are simple. Moreover, one can easily check that if 
$r\le 2^d$ then
$$x_{i+1}^{(r)}-\left(x_{i}^{(r)}-t\right)>2t,\,1\le i\le p+1.$$
Since $Q=P_{2^d}$ the above algorithm shows that $Q$ may be constructed from 
$P$ by using a total of $s=(p+1)2^d$ simple nondegenerate contractions.
\end{proof}

\begin{definition}\label{discrep}
Let $P(x)=\prod_{i=1}^{n}(x-x_i)$ and $Q(x)=\prod_{i=1}^{n}(x-y_i)$ be two 
hyperbolic polynomials of degree $n\ge 1$ whose zeros are arranged in 
nondecreasing order, so that if $n\ge 2$ then $x_i\le x_{i+1}$ and $y_i\le 
y_{i+1}$ for $1\le i\le n-1$. The number
$$\de(P,Q):=\sharp \left\{i\in \{1,\ldots,n\}\mid x_i\neq y_i\right\}$$
is called the {\em discrepancy} between $P$ and $Q$.
\end{definition}

\begin{remark}
It is clear from Definition~\ref{discrep} that $P=Q$ if and only if 
$\de(P,Q)=0$. 
\end{remark}

\begin{proof}[Proof of Proposition~\ref{contr1}]
The proposition is clearly true if $n=2$ and we may therefore assume that 
$n\ge 3$. Let $x_1<x_2<\ldots<x_n$ and $y_1<y_2<\ldots<y_n$ denote the 
zeros of $P$ and $Q$, respectively. Let further $r=\de(P,Q)$ and note that 
$r\ge 1$ since $P$ and $Q$ are distinct polynomials. Actually, since the 
condition $Q\pc P$ implies that 
$\sum_{i=1}^{n}y_i=\sum_{i=1}^{n}x_i$ we see that $r\ge 2$. We now prove 
the proposition by induction on $r$. If $r=2$ then by Theorem~\ref{crit} 
there exist indices $1\le i<j\le n$ such that $y_k=x_k$ whenever $k\neq i,j$ 
and $y_i=x_i+\si$ while $y_j=x_j-\si$ for some $\si\in \bR$ that satisfies
$$0<\si<\min\left(x_{i+1}-x_i,x_j-x_{j-1},\frac{x_j-x_i}{2}\right).$$
This means that if $j=i+1$ then $Q$ is already a simple nondegenerate 
contraction of $P$. If this is not the case then Lemma~\ref{c1} implies that 
$Q$ may be obtained from $P$ by the successive application 
of a finite number of simple nondegenerate contractions, which proves the 
result for $r=2$.

Suppose that $r\ge 3$ and assume that the proposition is true for all pairs 
of strictly hyperbolic polynomials whose discrepancies are at most $r-1$. 
Since $\sum_{i=1}^{n}(x_i-y_i)=0$ there must exist both positive and negative 
numbers among the differences $x_i-y_i$, $1\le i\le n$. A close examination 
of consecutive differences shows that at least one the following cases has to 
occur:

\smallskip

{\em Case 1.} There exists $i\in \{1,2,\ldots,n\}$ such that $x_i<y_i$ and 
$x_{i+1}>y_{i+1}$. Define the polynomial $R=\calT(i,i+1;t)P\in \calH_n$, where
$t=\min(y_i-x_i,x_{i+1}-y_{i+1})$. Note that 
$t\in \left(0,\frac{x_{i+1}-x_i}{2}\right)$ and thus $R$ is a simple 
nondegenerate contraction of $P$. We now use 
Theorem~\ref{crit} to check that one also has $Q\pc R$. This is obvious if 
$i=1$ and we may therefore assume that $i\ge 2$. It is then clear that
$$\sum_{k=1}^{m}x_k\le \sum_{k=1}^{m}y_k\text{ if }m\le i-1\text{ and }
\sum_{k=m}^{n}x_k\ge \sum_{k=m}^{n}y_k\text{ if }m\ge i+2.$$
Moreover, using the fact that $Q\pc P$ we get
\begin{equation*}
\begin{split}
& (x_i+t)+\sum_{k=1}^{i-1}x_k=y_i+\sum_{k=1}^{i-1}x_k\le \sum_{k=1}^{i}y_k
\text{ and }\\
& (x_{i+1}-t)+(x_i+t)+\sum_{k=1}^{i-1}x_k=\sum_{k=1}^{i+1}x_k\le 
\sum_{k=1}^{i+1}y_k,
\end{split}
\end{equation*}
which shows that if $t=y_i-x_i$ then the zeros of $Q$ and $R$ satisfy the 
inequalities in Theorem~\ref{crit}. It follows that $R$ is a strictly 
hyperbolic polynomials that satisfies $Q\pc R$ and $\de(Q,R)\le r-1$. Similar 
computations show that these relations remain true if $t=x_{i+1}-y_{i+1}$. 
By assumption, $Q$ may be obtained from $R$ by the successive 
application of a finite number of simple nondegenerate contractions. Since 
$R$ itself is a simple nondegenerate contraction of $P$, this proves the 
proposition in this case.

\smallskip

{\em Case 2.} There exist indices $i,j\in \{1,2,\ldots,n\}$ with $j\ge i+2$ 
such that $x_i<y_i$, $x_j>y_j$ and $x_k=y_k$ for $i+1\le k\le j-1$. Let 
$\si=\min(y_i-x_i,x_j-y_j)$ and set
$$R(x):=(x-x_i-\si)(x-x_j+\si)\prod_{\substack{k=1 \\ k\neq i,j}}^{n}(x-x_k),$$
so that $R$ is a strictly hyperbolic polynomial that satisfies $R\pc P$. 
Note that since $\si\in \left(0,\frac{x_j-x_i}{2}\right)$ it follows from 
Lemma~\ref{c1} that $R$ may be constructed by applying to $P$ a finite number 
of simple nondegenerate contractions. Clearly, these contractions affect only 
the zeros of $P$ and its successive transforms that lie in the interval
$[x_i,x_j]$. Computations similar to 
those used in case 1 show that $Q\pc R$. Moreover, it is clear that 
$\de(Q,R)\le r-1$. Using again the induction assumption we deduce that $Q$ 
may be obtained from 
$R$ and therefore also from $P$ by the successive application of a finite 
number of simple nondegenerate contractions, which completes the proof.
\end{proof}

Before proceeding with the proof of Theorem~\ref{main1} let us point out that 
if the nondegeneracy condition is 
omitted then minor modifications of the above arguments yield an 
analog of Proposition~\ref{contr1} for polynomials with multiple zeros. 
This result will not be used in the sequel and so we state it without proof:

\begin{proposition}\label{contr2}
Let $P$ and $Q$ be distinct polynomials in $\calH_{n}$ 
that satisfy $Q\pc P$. There exists a finite sequence of hyperbolic
 polynomials 
$P_{1},\ldots,P_{m}\in \calH_{n}$ such that $P_{1}=P$, $P_{m}=Q$ and 
$P_{i+1}$ is a simple contraction of $P_{i}$ for $1\le i\le m-1$.\hfill 
$\square$
\end{proposition}

The following proposition is the main step in the proof of Theorem~\ref{main1}.

\begin{proposition}\label{simplecontr}
If $P$ and $Q$ are strictly hyperbolic polynomials in $\calH_n$ such that 
$Q$ is a simple nondegenerate contraction of $P$ then $Q-\la Q'\pc P-\la P'$ 
for any $\la\in \bR$.
\end{proposition}

For the proof of Proposition~\ref{simplecontr} we need several additional 
results.
Let us first fix the notation that we shall use throughout this proof. 

\begin{nota}\label{not-gen}
We start with a strictly hyperbolic polynomial $P\in \calH_n$ given by
$$P(x)=\prod_{i=1}^{n}(x-x_i)\text{ and }P'(x)=n\prod_{j=1}^{n-1}(x-w_j).$$
By Rolle's theorem we may label the zeros of $P$ and $P'$ so that
$$x_1<w_1<x_2<\ldots<x_{n-1}<w_{n-1}<x_n,$$
which we assume henceforth. In most of the arguments below we shall also 
tacitly assume that $n\ge 3$. Fix an 
index $i\in \{1,2,\ldots,n-1\}$ and set 
$I=\left(0,\frac{x_{i+1}-x_i}{2}\right)$. For $t\in \bar{I}$ we let 
$P_t\in \calH_n$ denote the polynomial
$$P_t(x)
=(x-x_i-t)(x-x_{i+1}+t)\!\!\!\prod_{\substack{k=1 \\ k\neq i,i+1}}^{n}
\!\!\!(x-x_k)$$
and we define the following homotopy of polynomial pencils:
$$P(\la,t;x)=P_t(x)-\la P'_t(x), \text{ where }(\la,t)\in \bR\times \bar{I}
\text{ and }P'_t(x)=\frac{\pa}{\pa x}P_t(x).$$
Note that $P_t$ is a strictly hyperbolic polynomial whenever $t\in \{0\}\cup I$
and so by the Hermite-Poulain-Jensen theorem \cite[Theorem 5.4.9]{RS} 
the polynomial $P(\la,t;x)$ is strictly 
hyperbolic for all $(\la,t)\in \bR\times \left(\{0\}\cup I\right)$. Actually, 
if $0<\ve<\min(x_i-x_{i-1}, x_{i+2}-x_{i+1})$ then the same arguments 
show that the polynomial $P(\la,t;x)$ has only simple (real) zeros for any
$(\la,t)\in \bR\times \left(-\ve,\frac{x_{i+1}-x_i}{2}\right)$. If we now fix 
such an $\ve$ it follows from the implicit function theorem that the 
zeros of $P(\la,t;x)$ are real analytic functions of $(\la,t)$ in the domain 
$\bR\times \left(-\ve,\frac{x_{i+1}-x_i}{2}\right)$. Therefore, if we write
$$P(\la,t;x)=\prod_{k=1}^{n}\left(x-x_k(\la,t)\right)\!\text{ and }\!
P'(\la,t;x):=\frac{\pa}{\pa x}P(\la,t;x)=n\prod_{l=1}^{n-1}
\left(x-w_l(\la,t)\right)$$
and further assume that the zeros and the critical points of $P(\la,t;x)$ are 
labeled so 
that $x_k(0,0)=x_k$, $1\le k\le n$, and $w_l(0,0)=w_l$, $1\le l\le n-1$, then 
one has
\begin{equation}\label{int}
x_1(\la,t)<w_1(\la,t)<x_2(\la,t)<\ldots<x_{n-1}(\la,t)<w_{n-1}(\la,t)
<x_n(\la,t)
\end{equation}
if $(\la,t)\in \bR\times \left(\{0\}\cup I\right)$. These notations
will be used in all lemmas below.
\end{nota}

\begin{lemma}\label{l1}
If $1\le k\le n$ and $(\la,t)\in \bR\times \left(\{0\}\cup I\right)$ then 
$P'(\la,t;x_k(\la,t))\neq 0$ and
$$\frac{\pa}{\pa \la}x_k(\la,t)=\frac{P'_t(x_k(\la,t))}{P'(\la,t;x_k(\la,t))}
>0.$$
In particular, for all $j\in \{1,2,\ldots,n-1\}$ one has 
$$x_j(\la,t)<w_j(0,t)<x_{j+1}(\la,t)\text{ and }
\lim_{\la\rightarrow \infty}x_j(\la,t)
=\lim_{\la\rightarrow -\infty}\!x_{j+1}(\la,t)=w_j(0,t).$$ 
Moreover, $\lim_{\la\rightarrow \infty}x_n(\la,t)=
-\lim_{\la\rightarrow -\infty}x_1(\la,t)=\infty$.
\end{lemma}

\begin{proof}
The first assertion follows from the fact that $P(\la,t;x)$ is strictly 
hyperbolic and $P(\la,t;x_k(\la,t))=0$. Implicit differentiation with respect 
to $\la$ in the identity
$$P_t(x_k(\la,t))-\la P'_t(x_k(\la,t))=0$$
yields immediately the equality stated in the lemma. Note that since $P_t$ is 
strictly hyperbolic we have $P'_t(x_k(\la,t))\neq 0$, so that if we let
$P''_t(x)=\frac{\pa}{\pa x}P'_t(x)$ then
$$\left[P'_t(x_k(\la,t))\right]^2
\left[\frac{\pa}{\pa \la}x_k(\la,t)\right]^{-1}
\!=\left[P'_t(x_k(\la,t))\right]^2-P_t(x_k(\la,t))P''_t(x_k(\la,t))>0$$
by Laguerre's inequality for (strictly) hyperbolic polynomials 
\cite[Lemma 5.4.4]{RS}. If $t\in \{0\}\cup I$ is fixed then 
$-\la^{-1}P(\la,t;x)\rightarrow P'_t(x)$  as $|\la|\rightarrow \infty$ 
uniformly on compact 
sets. It follows that for $1\le j\le n-1$ one has 
$x_j(\la,t)<\lim_{\mu\rightarrow \infty}x_j(\mu,t)=w_j(0,t)$ and 
$x_{j+1}(\la,t)>\lim_{\mu\rightarrow -\infty}x_{j+1}(\mu,t)=w_j(0,t)$, which 
finishes the proof.
\end{proof}

For $1\le k\le n$ and $(\la,t)\in \bR\times \left(\{0\}\cup I\right)$ we 
define the following expressions: 
\begin{equation}\label{expres}
\begin{split}
F_k(\la,t) & =\left[\frac{P_t(x_k(\la,t))}
{\left(x_k(\la,t)-x_i-t\right)(x_k(\la,t)-x_{i+1}+t)
P'_t(x_k(\la,t))}\right]^2\text{ if } \la\neq 0,\\
F_i(0,t) & =F_{i+1}(0,t)=\frac{1}{(2t+x_i-x_{i+1})^2}\text{ and }
F_k(0,t)=0\text{ if }k\neq i,i+1.
\end{split}
\end{equation}
Note that $F_k(0,t)=\lim_{\la\rightarrow 0}F_k(\la,t)$ for all 
$k\in \{1,2,\ldots,n\}$ and $t\in \{0\}\cup I$.

\begin{lemma}\label{l2}
If $1\le k\le n$ and $(\la,t)\in \bR\times \left(\{0\}\cup I\right)$ then
$$\frac{\pa}{\pa t}x_k(\la,t)=(2x_k(\la,t)-x_i-x_{i+1})(2t+x_i-x_{i+1})
F_k(\la,t)\frac{\pa}{\pa \la}x_k(\la,t),$$
where $F_k(\la,t)$ is as in~\eqref{expres}.
\end{lemma}

\begin{proof}
By Lemma~\ref{l1} one has $\frac{\pa}{\pa \la}x_k(\la,t)\big|_{(0,t)}=1$ for 
all $t\in \{0\}\cup I$ and $1\le k\le n$. Moreover, it is clear that
$\frac{\pa}{\pa t}x_k(\la,t)\big|_{(0,t)}=0$ if $k\neq i,i+1$ while 
$\frac{\pa}{\pa t}x_i(\la,t)\big|_{(0,t)}=
-\frac{\pa}{\pa t}x_{i+1}(\la,t)\big|_{(0,t)}=1$. Thus, if $\la=0$ then the 
lemma is a consequence of~\eqref{expres}.

Assume now that $\la\neq 0$, so that
$$\frac{1}{\la}=\frac{P'_t(x_k(\la,t))}{P_t(x_k(\la,t))}=
\frac{1}{x_k(\la,t)-x_i-t}+\frac{1}{x_k(\la,t)-x_{i+1}+t}
+\!\!\sum_{\substack{r=1 \\ r\neq i,i+1}}^{n}\frac{1}{x_k(\la,t)-x_r}.$$
Applying $\frac{\pa}{\pa t}$ to the relation 
$P_t(x_k(\la,t))-\la P'_t(x_k(\la,t))=0$ we get
\begin{eqnarray*}
&& \left[P'_t(x_k(\la,t))-\la P''_t(x_k(\la,t))\right]
\frac{\pa}{\pa t}x_k(\la,t)=
\frac{\pa}{\pa t}\left[-P_t(x)+\la P'_t(x)\right]\Big|_{x=x_k(\la,t)}\\
&& =(2t+x_i-x_{i+1})\left[\prod_{\substack{r=1 \\ r\neq i,i+1}}^{n}
\!\!\!\left(x_k(\la,t)-x_r\right)-\la\!\!\!
\sum_{\substack{r=1 \\ r\neq i,i+1}}^{n}
\,\prod_{\substack{s=1 \\ s\neq i,i+1,r}}^{n}
\!\!\!\!\!\left(x_k(\la,t)-x_s\right)\right]\\
&& =(2t+x_i-x_{i+1})\left[1-\la\!\!\!\sum_{\substack{r=1 \\ r\neq i,i+1}}^{n}
\frac{1}{x_k(\la,t)-x_r}\right]\prod_{\substack{r=1 \\ r\neq i,i+1}}^{n}
\!\!\!\left(x_k(\la,t)-x_r\right)\\
&& =\frac{\la (2t+x_i-x_{i+1})P_t(x_k(\la,t))}
{(x_k(\la,t)-x_i-t)(x_k(\la,t)-x_{i+1}+t)}
\left[\frac{1}{\la}
-\sum_{\substack{r=1 \\ r\neq i,i+1}}^{n}\frac{1}{x_k(\la,t)-x_r}\right]\\
&& =\frac{(2x_k(\la,t)-x_i-x_{i+1})(2t+x_i-x_{i+1})P_t(x_k(\la,t))^2}
{(x_k(\la,t)-x_i-t)^2(x_k(\la,t)-x_{i+1}+t)^2P'_t(x_k(\la,t))}\\
&& =(2x_k(\la,t)-x_i-x_{i+1})(2t+x_i-x_{i+1})F_k(\la,t)P'_t(x_k(\la,t)).
\end{eqnarray*}
The result follows readily from Lemma~\ref{l1} since 
$P'_t(x_k(\la,t))\neq \la P''_t(x_k(\la,t)$.
\end{proof}

\begin{lemma}\label{l3}
Let $m\in\{1,2,\ldots,n\}$ and 
$(\la,t)\in \bR\times \left(\{0\}\cup I\right)$. Then
$$\sum_{k=1}^{m}x_k(\la,t)\ge \sum_{k=1}^{m}x_k(\la,0)\text{ if }m\le i-1,
\sum_{k=m}^{n}x_k(\la,t)\le \sum_{k=m}^{n}x_k(\la,0)\text{ if }m\ge i+2.$$
\end{lemma}

\begin{proof}
If $(\la,t)\in \bR\times \left(\{0\}\cup I\right)$ then~\eqref{int} and 
Lemma~\ref{l1} imply that
$$x_k(\la,t)<w_k(0,t)<x_{k+1}(0,t)\le x_i(0,t)<\frac{x_i+x_{i+1}}{2}$$
whenever $k\le i-1$ while for $k\ge i+2$ one gets that
$$x_k(\la,t)>w_{k-1}(0,t)>x_{k-1}(0,t)\ge x_{i+1}(0,t)>\frac{x_i+x_{i+1}}{2}.$$
Furthermore, by Lemma~\ref{l1} one has that 
$\frac{\pa}{\pa \la}x_k(\la,t)>0$ and by~\eqref{expres} we know that 
$F_k(\la,t)>0$ if $k\neq i,i+1$. Therefore, the above inequalities together 
with Lemma~\ref{l2} yield
$$\frac{\pa}{\pa t}x_k(\la,t)>0\text{ if } k\le i-1\text{ and }
\frac{\pa}{\pa t}x_k(\la,t)<0\text{ if } k\ge i+2.$$
It follows that all the inequalities in the lemma are strict if 
$(\la,t)\in \bR\times I$.
\end{proof}

We can now give a proof of Proposition~\ref{simplecontr}:

\begin{proof}[Proof of Proposition~\ref{simplecontr}]
Using the above notations we let $i\in \{1,2,\ldots,n\}$ and $\si\in I$ be 
such that $Q=\calT(i,i+1;\si)P$, so that
$$P(\la,0;x)=P(x)-\la P'(x)\text{ and } P(\la,\si;x)=Q(x)-\la Q'(x).$$
It is clear that for any $\la\in \bR$ one has
\begin{equation}\label{sum1}
\sum_{k=1}^{n}x_k(\la,0)=\sum_{k=1}^{n}x_k(\la,\si)=\sum_{k=1}^{n}x_k+n\la,
\end{equation}
where $x_k$, $1\le k\le n$, denote as before the zeros of $P$. By 
Theorem~\ref{crit} and~\eqref{sum1} we see that the relation 
$Q-\la Q'\pc P-\la P'$ is equivalent to
\begin{equation}\label{sum2}
\sum_{k=1}^{m}x_k(\la,0)\le \sum_{k=1}^{m}x_k(\la,\si),\,1\le m\le n-1.
\end{equation}
These inequalities are trivially true if $\la=0$ and so we may assume that
$\la\neq 0$. We distinguish two cases:

\smallskip

{\em Case 1: $\la>0$.} By Lemma~\ref{l1} one has
$\frac{\pa}{\pa \la}x_k(\la,t)>0$. Thus, if $\la>0$ then
$$x_{i+1}(\la,t)>x_{i+1}(0,t)=x_{i+1}-t>\frac{x_i+x_{i+1}}{2}\text{ for }
t\in [0,\si].$$
It follows from Lemma~\ref{l2} that $\frac{\pa}{\pa t}x_k(\la,t)<0$ if 
$\la >0$ and $t\in [0,\si]$. In particular,
\begin{equation}\label{case1}
x_{i+1}(\la,\si)<x_{i+1}(\la,0)\text{ if }\la>0.
\end{equation}

\smallskip

{\em Case 2: $\la<0$.} From Lemma~\ref{l1} again we deduce that in this case 
one has
$$x_i(\la,t)<x_i(0,t)=x_i+t<\frac{x_i+x_{i+1}}{2}\text{ for }
t\in [0,\si],$$
so that by Lemma~\ref{l2} we get $\frac{\pa}{\pa t}x_k(\la,t)>0$ if 
$\la<0$ and $t\in [0,\si]$. Hence
\begin{equation}\label{case2}
x_i(\la,\si)>x_i(\la,0)\text{ if }\la<0.
\end{equation}

Combining Lemma~\ref{l3} with~\eqref{case1} and~\eqref{case2} we see that for 
any $\la\in \bR\setminus \{0\}$ one has either
\begin{equation*}
\begin{split}
& \sum_{k=1}^{m}x_k(\la,0)\le \sum_{k=1}^{m}x_k(\la,\si),\,m\le i,
\sum_{k=m}^{n}x_k(\la,0)\ge \sum_{k=m}^{n}x_k(\la,\si),\,
m\ge i+2;\text{ or}\\
& \sum_{k=1}^{m}x_k(\la,0)\le \sum_{k=1}^{m}x_k(\la,\si),\,m\le i-1,
\sum_{k=m}^{n}x_k(\la,0)\ge \sum_{k=m}^{n}x_k(\la,\si),\,
m\ge i+1.
\end{split}
\end{equation*}
It is not difficult to see that these relations together with~\eqref{sum1} 
yield the inequalities in~\eqref{sum2}, which completes the proof of the 
proposition.
\end{proof}

Theorem~\ref{main1} is now an almost immediate consequence of the above 
results: 

\begin{proof}[Proof of Theorem~\ref{main1}]
In the generic case when both $P$ and $Q$ are strictly hyperbolic polynomials 
it follows 
from Proposition~\ref{contr1} that $Q$ may be obtained from $P$ by the 
successive application of a finite number of simple nondegenerate 
contractions. Therefore, in this case the theorem follows directly from 
Proposition~\ref{simplecontr}.

For the general case we let 
$x_1\le x_2\le\ldots\le x_n$ and $y_1\le y_2\le\ldots\le y_n$ denote the 
zeros of $P$ 
and $Q$, respectively, counted according to their respective multiplicities. 
Choose an arbitrary positive number $\ve$ and let $P_{\ve}$ 
and $Q_{\ve}$ be the polynomials with zeros $x_{i}-(n-i)\ve$, 
$1\le i\le n-1$, $x_{n}+\frac{n(n-1)}{2}\ve$, and $y_{i}-(n-i)\ve$, 
$1\le i\le n-1$, $y_{n}+\frac{n(n-1)}{2}\ve$, respectively. Note that both
$P_{\ve}$ and $Q_{\ve}$ are strictly hyperbolic and that 
$Q_{\ve}\pc P_{\ve}$. The above arguments imply that 
\begin{equation}\label{cont}
n^{-1}Q'_{\ve}\pc n^{-1}P'_{\ve}\text{ in }\calH_{n-1}\text{ and }
Q_{\ve}+\la Q_{\ve}'\pc P_{\ve}+\la P_{\ve}',\,\la \in \bR.
\end{equation}
Clearly, the zeros and the critical points of $P_{\ve}$ and $Q_{\ve}$ are 
continuous functions of $\ve$. The desired conclusion follows from 
Theorem~\ref{crit} and~\eqref{cont} by letting $\ve \rightarrow 0$.
\end{proof} 

\subsection{Applications to differential operators of Laguerre-P\'olya 
type and Appell polynomials}\label{ss22}

Theorem~\ref{main1} has several interesting consequences.
In order to state these we need some additional notations and definitions.

\begin{nota}\label{no1}
Given a nonconstant polynomial $P\in\vP$ we denote the 
barycenter of its zeros by $\fm(P)$. Suppose that
$$f(x)=\sum_{k=0}^{\infty}a_kx^k=x^mg(x),\quad x\in \bC,$$
is an entire function, where $m$ is a nonnegative integer and $g$ is an 
entire function such that $g(0)\neq 0$. One has a well-defined operator 
$f(D)\in \E\vP$ given by
$$f(D)[P](x)=\sum_{k=0}^{\infty}a_kP^{(k)}(x),\quad P\in \vP,$$
since only finitely many terms in this series are nonzero and so the lack of 
growth control on the coefficients in the power series expansion of $f$ 
causes no problems. We associate to $f$ an infinite family of differential 
operators $\left\{\calD(f,n)\right\}_{n=m+1}^{\infty}$ defined as follows:
\begin{equation}\label{inffam}
\calD(f,n)=k_n(f)f(D),\text{ where } 
k_n(f)=\left[\binom{n}{m}f^{(m)}(0)\right]^{-1},\, n\ge m+1.
\end{equation}
Note that these operators are in fact rescalings of $f(D)$ chosen so that if 
$n\ge m+1$ then $\calD(f,n)$ maps monic polynomials of degree 
$n$ to monic polynomials of degree $n-m$. In particular, if $m=0$ then all 
operators $\calD(f,n)$, $n\in \bN$, coincide with $f(0)^{-1}\!f(D)$ and 
preserve the class of monic polynomials of degree $d$ for any $d\ge 0$. 
\end{nota}  

\begin{definition}\label{LP}
A real entire function $\vf$ is said to be in the {\em Laguerre-P\'olya class},
$\vf\in \calLP$, if it can be expressed in the form
\begin{equation}\label{genLP}
\vf(x)=cx^me^{-a^2x^2+bx}\prod_{k=1}^{\infty}(1-\al_k x)e^{\al_k x},\quad 
x\in \bC,
\end{equation}
where $a,b,c,\al_k\in \bR$, $c\neq 0$, $m$ is a nonnegative integer,
$\sum_{k=1}^{\infty}\al_{k}^{2}<\infty$ and where, by the usual convention, 
the canonical product reduces to $1$ if $\al_k=0$ for all $k\in \bN$. An 
operator $T\in \E\vP$ is said to be a {\em differential operator of 
Laguerre-P\'olya type} if $T=\vf(D)$, where $\vf\in \calLP$. 
\end{definition}

\begin{nota}\label{not2}
Let $\calLP_0:=\left\{\vf \in \calLP\mid \vf(0)\neq 0\right\}$. For 
$m\in \bN$ we set
\begin{equation*}
\calLP_m=x^m\calLP_0=\left\{\vf \in \calLP\mid \vf^{(k)}(0)=0,\,0\le k\le m-1,
\,\vf^{(m)}(0)\neq 0\right\}.
\end{equation*}
Clearly, $\calLP$ is a commutative monoid under ordinary multiplication of 
functions. Actually, $\calLP$ may be viewed as a $\bZ_+$-graded monoid, where 
$\bZ_+$ denotes the additive monoid of nonnegative integers. Indeed, note 
that $\calLP_0$ is a submonoid of $\calLP$ which acts on $\calLP_m$ for each 
$m\in \bZ_+$ and that $\calLP$ decomposes into a disjoint union
\begin{equation}\label{disj}
\calLP=\bigcup_{m=0}^{\infty}\calLP_m\text{ with } 
\calLP_{m_1}\cdot\calLP_{m_2}=\calLP_{m_1+m_2}\text{ for }m_1,m_2\in \bZ_+.
\end{equation}
\end{nota}

As we already pointed out in the introduction, by a classical 
theorem of P\'olya one knows that all differential operators of 
Laguerre-P\'olya type map hyperbolic polynomials to
hyperbolic polynomials. By using Theorem~\ref{main1} one can actually 
show that all such operators are in fact natural preservers of the 
spectral order:

\begin{corollary}\label{iso}
Let $m,n\in \bZ_+$ with $n\ge m+1$ and $\vf\in \calLP_m$. If $P,Q\in \calH_n$ 
are such that $Q\pc P$ then $\calD(\vf,n)[Q]\pc \calD(\vf,n)[P]$ in 
$\calH_{n-m}$.
\end{corollary}

\begin{remark}\label{nontriv}
It is clear that if $\vf\in \calLP_m$ then $\calD(\vf,m)[P]\equiv 1$ for all 
$P\in \calH_m$ while $\calD(\vf,n)[P]\equiv 0$ if $P\in \calH_n$ with 
$n\le m-1$. This is the reason why we impose the condition $n\ge m+1$ both in 
Corollary~\ref{iso} and Corollary~\ref{monot} of \S \ref{s3}.
\end{remark}

To prove Corollary~\ref{iso} we need to establish first the following result.

\begin{lemma}\label{deriv}
Let $n\ge 2$ and $P,Q\in \calH_n$ with $Q\pc P$. Then 
$n^{-1}Q'\pc n^{-1}P'$ in $\calH_{n-1}$.
\end{lemma}

\begin{proof}
It is enough to prove the lemma in the generic case when 
$P$ and $Q$ are strictly hyperbolic polynomials and $Q$ is a simple 
nondegenerate contraction of $P$ (the general case follows from this one 
by arguing as in the proof of Theorem~\ref{main1}). Let then 
$Q=\calT(i,i+1;\si)P$, where $\si\in I$ and $i\in \{1,2,\ldots,n\}$.
Using Notation~\ref{not-gen} we may write
\begin{equation*}
\begin{split}
& P(\la,0;x)=P(x)-\la P'(x)=\prod_{k=1}^{n}(x-x_k(\la,0)),\,\,
P'(\la,0;x)=n\prod_{l=1}^{n-1}(x-w_l(\la,0)),\\
& P(\la,\si;x)=Q(x)-\la Q'(x)=\prod_{k=1}^{n}(x-x_k(\la,\si)),\,\,
P'(\la,\si;x)=n\prod_{l=1}^{n-1}(x-w_l(\la,\si)).
\end{split}
\end{equation*}
By Proposition~\ref{simplecontr} we know that $P(\la,\si;x)\pc P(\la,0;x)$, so 
that~\eqref{sum2} is valid. Therefore, if we let $\la\rightarrow \infty$ 
in~\eqref{sum2} and use the second part of Lemma~\ref{l1} we obtain
\begin{equation}\label{der}
\sum_{j=1}^{m}w_j(0,0)\le \sum_{j=1}^{m}w_j(0,\si),\,1\le m\le n-1.
\end{equation} 
Since $Q$ is a contraction of $P$ one has 
$Q\pc P$, so that $\fm(Q)=\fm(P)$ and thus $\fm(Q')=\fm(P')$. This shows that 
the inequality in~\eqref{der} corresponding to $m=n-1$ is actually an equality,
which by Theorem~\ref{crit} proves the lemma.
\end{proof}

\begin{proof}[Proof of Corollary~\ref{iso}]
Let $X=(x_1,x_2,\ldots,x_n)$ and $Y=(y_1,y_2,\ldots,y_n)$ be two unordered 
$n$-tuples of real numbers and set
\begin{equation*}
d(X,Y)=\min_{\pi\in \Sigma_n}\max_{1\le i\le n}\left|x_i-y_{\pi(i)}\right|.
\end{equation*}
This is the so-called {\em optimal matching distance} between the unordered 
$n$-tuples $X$ and $Y$. It is not difficult to see that $d$ defines a metric
 on the quotient space 
$\bR^n/\Sigma_n$ of all such $n$-tuples and therefore also on the manifold 
$\calH_n$ in view of~\eqref{rmap}.

We use the rearrangement-free 
characterization of the spectral order given in Theorem~\ref{spec} (i) in the 
following way: to any function $f:\bR\rightarrow \bR$ we associate a 
function $\tilde{f}:\bR^n/\Sigma_n\rightarrow \bR$ by setting 
\begin{equation}\label{tild}
\tilde{f}(X)=\sum_{i=1}^{n}f(x_i)\text{ for }
X=(x_1,x_2,\ldots,x_n)\in \bR^n/\Sigma_n.
\end{equation} 
If $f$ is convex then Theorem~\ref{spec} (i) asserts that 
$\tilde{f}(X)\le \tilde{f}(Y)$ whenever $X\pr Y$, that is, $\tilde{f}$ is a 
{\em Schur-convex} function (cf.~\cite[Ch.~3]{MO}). Thus 
$X\pr Y$ if and only if $\tilde{f}(X)\le \tilde{f}(Y)$ 
for any function $\tilde{f}$ as in~\eqref{tild} associated to a convex 
function $f$. 

Assume now that 
$P,Q\in \calH_n$ are such that $Q\pc P$ and let $\vf\in \calLP_m$, where 
$m\in \bZ_+$, $m\le n-1$. Suppose that $\vf$ is as in~\eqref{genLP}
with Maclaurin expansion 
$$\vf(x)=\sum_{k=m}^{\infty}a_kx^k,\quad x\in \bC.$$
For $j\in \bN$ let $\tau_j:=b+\sum_{\nu=1}^{j}\al_\nu$ and define the 
following polynomials:
\begin{equation}\label{approx}
\vf_j(x)=cx^m\left(1-\frac{ax}{\sqrt{j}}\right)^{j}
\left(1+\frac{ax}{\sqrt{j}}\right)^{j}
\left(1+\frac{\tau_{j}x}{n_j}\right)^{n_j}\prod_{\nu=1}^{j}(1-\al_{\nu}x).
\end{equation}
It is a well-known fact that if one chooses $\{n_j\}_{j\in\bN}$ as a sequence 
of integers growing sufficiently fast to infinity as $j\rightarrow \infty$ 
then the sequence of hyperbolic polynomials $\{\vf_j\}_{j\in\bN}$ satisfies 
$\vf_j\rightrightarrows \vf$ as $j\rightarrow \infty$, where 
$\rightrightarrows$ denotes uniform 
convergence on all compact subsets of $\bC$ (see, e.g., \cite[Ch.~8]{L}). 
Therefore, if we let $N_j:=\deg\vf_j$ and write the polynomial $\vf_j$ as
$$\vf_j(x)=\sum_{k=m}^{\infty}a_{j,k}x^k,\quad x\in \bC,$$
with $a_{j,k}=0$ for $k\ge N_j+1$ then it follows from Cauchy's integral 
formula that $\lim_{j\rightarrow\infty}a_{j,k}=a_k$ for all $k\ge m$. This 
implies that for any fixed polynomial $R\in \vP$ with $\deg R=n$ one has
$$\vf_j(D)[R]=\sum_{k=m}^{n}a_{j,k}R^{(k)}\rightrightarrows 
\sum_{k=m}^{n}a_kR^{(k)}=\vf(D)[R]\text{ as }j\rightarrow \infty.$$
In particular, 
$\calD(\vf_j,n)[P]\rightrightarrows \calD(\vf,n)[P]$ and 
$\calD(\vf_j,n)[Q]\rightrightarrows \calD(\vf,n)[Q]$ as 
$j\rightarrow \infty$, so that
\begin{equation}\label{lim1}
\begin{split}
& d\big(\calZ(\calD(\vf_j,n)[P]),\calZ(\calD(\vf,n)[P])\big)
\longrightarrow 0\text{ and }\\
& d\big(\calZ(\calD(\vf_j,n)[Q]),\calZ(\calD(\vf,n)[Q])\big)
\longrightarrow 0\text{ as }j\longrightarrow\infty.
\end{split}
\end{equation}
On the other hand, by Theorem~\ref{main1} and Lemma~\ref{deriv} we know that
$$\calZ(\calD(\vf_j,n)[Q])\pr\calZ(\calD(\vf_j,n)[P])\text{ in }
\bR^{n-m}/\Sigma_{n-m}\text{ for }j\in \bN.$$
Thus, if $f$ is a real-valued convex function on $\bR$ and 
$\tilde{f}$ is as in~\eqref{tild} then 
\begin{equation}\label{lim2}
\tilde{f}\big(\calZ(\calD(\vf_j,n)[Q])\big)\le 
\tilde{f}\big(\calZ(\calD(\vf_j,n)[P])\big)\text{ for }j\in \bN.
\end{equation}
Since $f$ is convex on $\bR$ it is also continuous there 
and so $\tilde{f}$ is a continuous function on $\bR^n/\Sigma_n$. Therefore, 
by letting $j\rightarrow\infty$ in~\eqref{lim1} and~\eqref{lim2} we obtain
$$\tilde{f}\big(\calZ(\calD(\vf,n)[Q])\big)\le 
\tilde{f}\big(\calZ(\calD(\vf,n)[P])\big).$$
As explained above, this implies that 
$$\calZ(\calD(\vf,n)[Q])\pr\calZ(\calD(\vf,n)[P])\text{ in }
\bR^{n-m}/\Sigma_{n-m}.$$
Hence $\calD(\vf,n)[Q]\pc \calD(\vf,n)[P]$ in $\calH_{n-m}$, which
completes the proof.
\end{proof}

\begin{nota}\label{not3}
Define the following monoids of linear operators:
\begin{equation}\label{mono1}
\calA=\bigcap_{n=0}^{\infty}\calA_n,\text{ where }
\calA_n=\left\{T\in \E\vP\mid T\big(\calH_n\big)\subseteq \calH_n\right\},\,
n\in \bZ_+.
\end{equation}
Note that $\calA_n$ is the largest submonoid of $\E\vP$ consisting of 
linear operators that act on 
$\calH_n$ for fixed $n\in \bZ_+$, while $\calA$ is the largest submonoid of 
$\E\vP$ acting on each of the manifolds $\calH_n$, $n\in \bZ_+$. 
\end{nota}

In \cite[Theorem~1]{CPP} it was shown that 
\begin{equation}\label{mono2}
\calA=\left\{\vf(D)\mid \vf\in \calLP,\,\vf(0)=1\right\}\subset \calLP_0.
\end{equation}
From Corollary~\ref{iso} and~\eqref{mono2} we deduce that all operators in 
$\calA$ are {\em isotonic} (see Definition~\ref{d-set} below)
with respect to the spectral order on $\calH_n$ for any $n\in \bN$:

\begin{corollary}\label{coroiso1}
If $n\ge 1$ and $P,Q\in \calH_n$ are such that $Q\pc P$ then $T[Q]\pc T[P]$ 
for all operators $T\in \calA$.\hfill$\square$
\end{corollary}

Yet another consequence of Corollary~\ref{iso} is that 
the sequence of nonconstant Appell polynomials associated to any given function
in the Laguerre-P\'olya class may be characterized by means of a global 
minimum property 
with respect to the spectral order. Indeed, let $n\in \bN$ and consider the 
following submanifold of $\calH_n$:
\begin{equation}\label{h0}
\calH_{n}^{0}=\{P\in \calH_{n}\mid \fm(P)=0\}.
\end{equation}
Given $\vf\in \calLP$ and $n\in \bZ_+$ one defines the $n$-th {\em Appell 
polynomial} $g_{n}^{*}$ associated with $\vf$ by $g_{n}^{*}(x)=\vf(D)[x^n]$ 
(see, e.g., \cite{CC1}). Recall the decomposition of $\calLP$ 
from~\eqref{disj} and assume that $\vf\in \calLP_m$
for some $m\in \bZ_+$. Clearly, $g_{n}^{*}$ is a nonconstant polynomial if 
and only if $n\ge m+1$ (cf.~Remark~\ref{nontriv}). Corollary~\ref{iso}, 
Theorem~\ref{crit}, and the fact that 
$x^n\pc P(x)$ for any $P\in \calH_{n}^{0}$, $n\in \bN$, yield the following:

\begin{corollary}\label{coroiso2}
Let $m\in \bZ_+$ and $\vf\in \calLP_m$. If $n\ge m+1$ then the monic 
polynomial $k_n(\vf)g_{n}^{*}$ is the (unique) global minimum of the poset 
$\left(\calD(\vf,n)[\calH_{n}^{0}],\pc\right)$, where 
$\calD(\vf,n)[\calH_{n}^{0}]:=\{\calD(\vf,n)[P]\mid P\in \calH_{n}^{0}\}$, 
$k_n(\vf)$ is as in~\eqref{inffam} and $g_{n}^{*}$ is the $n$-th Appell 
polynomial associated with $\vf$.\hfill$\square$
\end{corollary}

In view of Theorems~\ref{spec} and~\ref{crit}, 
Corollary~\ref{coroiso2} admits the following 
geometrical interpretation: up to a factor $k_n(\vf)$ the $n$-th Appell 
polynomial associated with $\vf$ coincides with the (unique) 
polynomial in the image set
$\calD(\vf,n)[\calH_{n}^{0}]$ whose zeros are less spread out than the zeros 
of any other polynomial in this set.

\begin{remark}\label{topo}
A systematic 
investigation of the topological properties of $\calH_{n}$ 
and $\calH_{n}^{0}$ was initiated by Arnold in \cite{Ar}. These manifolds 
have since been extensively studied in singularity theory and related topics. 
\end{remark}

\section{Theorem~\ref{main2} and some consequences}\label{s3}

\subsection{Proof of Theorem~\ref{main2}}\label{ss31}

The result holds trivially for $n=1$ and so we may assume 
that $n\ge 2$. As in \S \ref{s2}, we start with a strictly hyperbolic 
polynomial $P\in\calH_n$ given by
$$P(x)=\prod_{i=1}^{n}(x-x_i)\text{ and }P'(x)=n\prod_{j=1}^{n-1}(x-w_j)$$
with $x_1<w_1<x_2<\ldots<x_{n-1}<w_{n-1}<x_n$ and we define the following 
pencils of polynomials:
$$P_{\la}(x)=P(x)-\la P'(x)\text{ and }P_{\la}'(x)=P'(x)-\la P''(x),\quad 
\la\in \bR.$$
Denote the zeros of $P_{\la}$ and $P_{\la}'$ by $x_i(\la)$, $1\le i\le n$, 
and $w_j(\la)$, $1\le j\le n-1$, respectively. If we assume that these 
are labeled so that $x_i(0)=x_i$, $1\le i\le n$, and $w_j(\la)=w_j$, 
$1\le j\le n-1$, then by letting $t=0$ in~\eqref{int} we see that
\begin{equation}\label{int2}
x_1(\la)<w_1(\la)<x_2(\la)<\ldots <x_{n-1}(\la)<w_{n-1}(\la)<x_n(\la)
\end{equation}
for all $\la\in \bR$. The following proposition is the key step in the proof 
of Theorem~\ref{main2}.

\begin{proposition}\label{incr}
If $P$ is as above then each of the functions $f_m:\bR\rightarrow \bR$ given by
$$f_m(\la)=\sum_{i=1}^{m}(x_i(\la)-\la),\quad 1\le m\le n-1,$$
is increasing on $(-\infty,0]$ and decreasing on $[0,\infty)$.
\end{proposition}

The proof of Proposition~\ref{incr} is based on two lemmas:

\begin{lemma}\label{lincr1}
Let $1\le j\le n-1$ and $\la\in \bR$. Then
$$\sum_{i=1}^{m}\frac{1}{x_i(\la)-w_j(\la)}<0$$
for all $m\in \{1,\ldots,n-1\}$.
\end{lemma}

\begin{proof}
If $m\le j$ then for each $i\le m$ one has $x_i(\la)\le x_m(\la)<w_j(\la)$ 
by~\eqref{int2}, so that in this case all terms in the sum are negative. 
Assume that $m\ge j+1$. Then
$$0=\frac{P_{\la}'(w_j(\la))}{P_{\la}(w_j(\la))}
=\sum_{i=1}^{n}\frac{1}{w_j(\la)-x_i(\la)}
=\sum_{i=1}^{m}\frac{1}{w_j(\la)-x_i(\la)}+
\sum_{i=m+1}^{n}\frac{1}{w_j(\la)-x_i(\la)}.$$
Thus
$$\sum_{i=1}^{m}\frac{1}{x_i(\la)-w_j(\la)}
=\sum_{i=m+1}^{n}\frac{1}{w_j(\la)-x_i(\la)}<0$$
since~\eqref{int2} implies that $x_i(\la)\ge x_{m+1}(\la)>w_j(\la)$ if 
$i\ge m+1$.
\end{proof}

\begin{lemma}\label{lincr2}
If $1\le j\le n-1$ and $\la\in \bR$ then
$$w_{j}'(\la)=\frac{P''(w_j(\la))}{P_{\la}''(w_j(\la))}>0,$$
where $P_{\la}''(x)=\frac{\pa}{\pa x}P_{\la}'(x)$.
\end{lemma}

\begin{proof}
Apply Lemma~\ref{l1} to $P'(\la,t,w_j(\la,t))$, $1\le j\le n-1$, 
and set $t=0$.
\end{proof}

\begin{proof}[Proof of Proposition~\ref{incr}]
Using Lemma~\ref{l1} and a partial fractional decomposition we get
$$x_{i}'(\la)-1=\frac{\la P''(x_i(\la))}{P_{\la}'(x_i(\la))}
=\sum_{j=1}^{n-1}\frac{P''(w_j(\la))}{P_{\la}''(w_j(\la))}
\frac{\la}{x_i(\la)-w_j(\la)}
=\sum_{j=1}^{n-1}\frac{\la w_{j}'(\la)}{x_i(\la)-w_j(\la)}.$$
Therefore, if $1\le m\le n-1$ then
\begin{equation}\label{ident}
f_{m}'(\la)=\sum_{i=1}^{m}(x_{i}'(\la)-1)=\la\sum_{j=1}^{n-1}
\sum_{i=1}^{m}\frac{w_{j}'(\la)}{x_i(\la)-w_j(\la)}.
\end{equation}
Lemmas~\ref{lincr1} and ~\ref{lincr2} imply that
$$\sum_{i=1}^{m}\frac{w_{j}'(\la)}{x_i(\la)-w_j(\la)}<0,\quad \la\in \bR,$$
which together with~\eqref{ident} shows that $\la f_{m}' (\la)<0$ if 
$\la\neq 0$, as required.
\end{proof}

Theorem~\ref{main2} is now a straightforward consequence of 
Theorem~\ref{crit} and the following result.

\begin{proposition}\label{allincr}
Let $P\in \calH_n$ and set $P_{\la}(x)=P(x)-\la P'(x)$, where $\la\in \bR$. 
For any fixed $\la$ denote the zeros of $P_{\la}$ by $x_i(\la)$, 
$1\le i\le n$, and arrange these so that $x_1(\la)\le\ldots\le x_n(\la)$. 
Given $m\in\{1,2,\ldots,n\}$ we define a function $f_m:\bR\rightarrow\bR$ by 
$$f_m(\la)=\sum_{i=1}^{m}(x_i(\la)-\la).$$
If $1\le m\le n-1$ then $f_m$ is nondecreasing on $(-\infty,0]$ and it is 
nonincreasing on $[0,\infty)$. Moreover, $f_n$ is a constant function on 
$\bR$.
\end{proposition}

\begin{proof}
The first assertion follows from Proposition~\ref{incr} since 
$P$ may be approximated by strictly hyperbolic polynomials in $\calH_n$ 
uniformly on compact subsets of $\bC$. Indeed, if  
$\ve\in\bR\setminus \{0\}$ then $\hat{P}_{\ve}(x):=(1-\ve D)^{n-1}P(x)$ is a 
strictly hyperbolic polynomial in $\calH_n$
(cf., e.g., \cite[Lemma 4.2]{CC2}). It is clear that 
$\hat{P}_{\ve}\rightrightarrows P$ as $\ve\rightarrow 0$. The second statement 
follows from the 
fact that $f_n(\la)=\sum_{i=1}^{n}x_i$ for all $\la\in \bR$, where $x_i$, 
$1\le i\le n$, are the zeros of $P$. 
\end{proof}

\begin{remark}\label{better}
Proposition~\ref{allincr} has recently been  
extended to arbitrary hyperbolic polynomial pencils in~\cite{BP}, where 
it was furthermore shown that $f_m$, $1\le m\le n-1$, are actually concave 
functions on $\bR$.
Note that by \cite[Theorem~4]{B} these partial sums cannot have a 
common local maximum unless the polynomial pencil under consideration is 
of logarithmic derivative type, i.e., of the form $P-\la P'$, $\la\in \bR$.
\end{remark}

\begin{corollary}\label{maincoro}
Let $\la_1,\la_2\in\bR$ be such that $\la_1\la_2\ge 0$ and 
$|\la_1|\le |\la_2|$. If $m,n\in\bZ_+$ with $n\ge\max(2,m+1)$ then for any 
$P\in \calH_n$ one has
$$\binom{n}{m}^{-1}D^{m}(1-\la_1 D)e^{\la_1D}P\pc 
\binom{n}{m}^{-1}D^{m}(1-\la_2 D)e^{\la_2D}P\text{ in }\calH_{n-m}.$$
In particular, if  $s_1,s_2\in\bR$ satisfy 
$s_1s_2\ge 0$ and $|s_1|\le |s_2|$ then
\begin{equation*}
\begin{split}
\binom{n}{m}^{-1}D^{m}(1-s_1\la D)e^{s_1\la D}P & \pc 
\binom{n}{m}^{-1}D^{m}(1-s_2\la D)e^{s_2\la D}P\\
e^{-s_{1}^2\la^2 D^2}P & \pc e^{-s_{2}^2\la^2 D^2}P
\end{split}
\end{equation*}
for all $P\in \calH_n$ and $\la\in \bR$.
\end{corollary}

\begin{proof}
The first relation is an immediate consequence of Theorem~\ref{crit}, 
Proposition~\ref{allincr} and repeated use of Lemma~\ref{deriv} since 
$(1-\la D)e^{\la D}P(x)=P(x+\la)-\la P'(x+\la)$ for all $\la\in \bR$.
Setting $\la_i=s_i\la$, $i=1,2$, one gets the second relation. Let 
$j\in \bN$ and define a function
$$\psi_j(x)=\left(1-\frac{\la^2 x^2}{j}\right)^j
=\left[\left(1-\frac{\la x}{\sqrt{j}}\right)
e^{^{\frac{\la x}{\sqrt{j}}}}\right]^j
\left[\left(1+\frac{\la x}{\sqrt{j}}\right)
e^{^{-\frac{\la x}{\sqrt{j}}}}\right]^j,$$
where $\la$ is a fixed real number. Clearly, the second relation implies that 
for any $P\in \calH_n$ and $j\in \bN$ one has 
$\psi_j(s_{1}D)[P]\pc \psi_j(s_{2}D)[P]$.
Moreover, from $\psi_j(x)\rightrightarrows e^{-\la^2 x^2}$ as 
$j\rightarrow\infty$ one easily gets $\psi_j(s_{i}D)[P]\rightrightarrows 
e^{-s_{i}^2\la^2 D^2}P$ for $i=1,2$. The third relation is obtained 
by letting $j\rightarrow\infty$.
\end{proof}

\subsection{Orbits of hyperbolic polynomials}\label{ss32}

Theorem~\ref{main2} and Corollary~\ref{maincoro} allow us to study the orbits  
of hyperbolic polynomials under the action of differential operators of 
Laguerre-P\'olya type. To do this we need some new notation.

\begin{nota}\label{not-orbit}
Let $\li$ denote the Banach algebra of bounded real sequences of the form 
$\{s_i\}_{i=0}^{\infty}$. We endow $\li$ with a partial ordering $\leqslant$ 
defined as follows: given two elements $\bs=\{s_i\}_{i=0}^{\infty}$ and 
$\bt=\{t_i\}_{i=0}^{\infty}$ of $\li$ we set 
$\bs\leqslant \bt$ if $|s_i|\le |t_i|$ and $s_it_i\ge 0$ for all 
$i\in \bZ_+$. For fixed $\bs=\{s_i\}_{i=0}^{\infty}\in \li$, $m\in \bZ_+$ and
a function $\vf\in \calLP_m$ of the form \eqref{genLP} we define the 
$\bs$-{\em deformation} of $\vf$ to be 
\begin{equation}\label{action}
\vf^{\bs}(x)=cx^me^{-s_{0}^2a^2x^2+bx}\prod_{k=1}^{\infty}(1-s_k\al_k x)
e^{s_k\al_k x},\quad x\in \bC.
\end{equation}
Note that $\vf^{\bs}\in \calLP_m$ and so \eqref{action} defines an action of 
$\li$ on $\calLP_m$ for any $m\in \bZ_+$
\begin{equation}\label{act}
\begin{split}
\li\times \calLP_m & \longrightarrow \calLP_m\\
(\bs,\vf) & \longmapsto \bs\cdot\vf:=\vf^{\bs}
\end{split}
\end{equation}
by means of which we associate to any $\vf\in \calLP_m$ an 
infinite-parameter family of deformations of the operator $\vf(D)$, namely  
$$\calF_{\vf}
:=\left\{\calD\left(\vf^{\bs},n\right)\mid \bs\in \li,\,n\in \bN,\,
n\ge m+1\right\},$$ 
where $\calD\left(\vf^{\bs},n\right)$ is as in \eqref{inffam}. 
\end{nota}

The operator families $\calF_{\vf}$ satisfy the following 
global monotony property with respect to the 
partial orderings $\leqslant$ on $\li$ and $\pc$ on $\calH_n$, $n\in \bZ_+$, 
respectively:

\begin{corollary}\label{monot}
Let $m,n\in \bZ_+$ with $n\ge m+1$ and $\vf\in \calLP_{m}$. If 
$\bs,\bt\in \li$ are such that $\bs\leqslant \bt$ then 
$\calD\left(\vf^{\bs},n\right)[P]\pc \calD\left(\vf^{\bt},n\right)[P]$ in 
$\calH_{n-m}$ for any $P\in \calH_n$.
\end{corollary}

\begin{proof}
Let us fix $\bs=\{s_i\}_{i=0}^{\infty}\in\li$ and 
$\bt=\{t_i\}_{i=0}^{\infty}\in\li$ such that 
$\bs\leqslant \bt$. Given $m,n\in\bZ_+$ with $n\ge\max(2,m+1)$ and 
$\vf\in\calLP_{m}$ as in~\eqref{genLP} we approximate 
$\vf^{\bs}(x)$ and $\vf^{\bt}(x)$ uniformly on compact subsets 
of $\bC$ by means of the functions
\begin{equation*}
\begin{split}
\vf_{j}^{\bs}(x) & =cx^{m}e^{-s_{0}^{2}a^{2}x^{2}+bx}
\prod_{k=1}^{j}(1-s_{k}\al_{k}x)e^{s_{k}\al_{k}x}\text{ and}\\ 
\vf_{j}^{\bt}(x) & =cx^{m}e^{-t_{0}^{2}a^{2}x^{2}+bx}
\prod_{k=1}^{j}(1-t_{k}\al_{k}x)e^{t_{k}\al_{k}x},
\end{split}
\end{equation*}
respectively, where $j\in \bN$. By Corollary~\ref{maincoro} we know that
\begin{equation}\label{last}
\calD(\vf_{j}^{\bs},n)[P]\pc \calD(\vf_{j}^{\bt},n)[P]\text{ in }
\calH_{n-m}
\end{equation}
for arbitrarily fixed $P\in \calH_n$ and $j\in \bN$. Standard arguments 
involving the uniform 
convergence of the above sequences of functions similar to those given in 
the proof of Corollary~\ref{iso} 
show that $\calD(\vf_{j}^{\bs},n)[P]\rightrightarrows 
\calD(\vf^{\bs},n)[P]$ and $\calD(\vf_{j}^{\bt},n)[P]\rightrightarrows 
\calD(\vf^{\bt},n)[P]$ as $j\rightarrow\infty$. The desired result 
follows from~\eqref{last} by letting $j\rightarrow\infty$.
\end{proof}

Recall from~\eqref{mono2} that $\calA$ is the largest submonoid of $\E\vP$ 
acting on each of the manifolds $\calH_n$, $n\in \bZ_+$. We define a 
binary relation on $\calA$ which by abuse of notation we denote again by $\pc$
in the following manner: given $T_1,T_2\in \calA$ set $T_1\pc T_2$ if
$T_{1}[P]\pc T_{2}[P]$ for all $P\in \calH_n$, $n\in \bN$.

\begin{lemma}\label{part-ord}
The pair $(\calA,\pc)$ is a poset.
\end{lemma}

\begin{proof}
Clearly, $\pc$ inherits the reflexivity and transitivity properties 
from the partial orderings on the posets $(\calH_n,\pc)$, $n\in \bZ_+$. 
Assume that $T_1,T_2\in \calA$ are such 
that $T_1\pc T_2$ and $T_2\pc T_1$. By~\eqref{mono1} we may write 
$T_i=\vf_i(D)$, where $\vf_i\in \calLP$ with $\vf_i(0)=1$, $i=1,2$. 
In particular, 
$\vf_1(D)[x^n]\pc \vf_2(D)[x^n]$ and $\vf_2(D)[x^n]\pc \vf_1(D)[x^n]$,
$n\in \bZ_+$. Since $(\calH_n,\pc)$ is a poset for all $n\in \bZ_+$ we deduce 
that the sequences of Appell polynomials associated to $\vf_1$ and $\vf_2$ 
must coincide. It follows that $\vf_1=\vf_2$ and thus $T_1=T_2$, which shows 
that $\pc$ is also antisymmetric.
\end{proof}

From Corollary~\ref{monot} we 
deduce the following compatibility relation between the posets 
$(\li,\leqslant)$ and $(\calA,\pc)$.

\begin{corollary}\label{coromonot3}
If $T\in \calA$ and $\bs,\bt\in \li$ with $\bs\leqslant \bt$ then 
$\bs\cdot T\pc \bt\cdot T$.\hfill$\square$
\end{corollary}

Let $\calLP'$ be the class of entire functions of the form
\begin{equation}\label{LP'}
\vf(x)=cx^me^{-a^2x^2}\prod_{k=1}^{\infty}(1-\al_k x)e^{\al_k x},\quad 
x\in \bC,
\end{equation}
where $a,c,\al_k\in \bR$, $c\neq 0$, $m\in \bZ_+$ and 
$\sum_{k=1}^{\infty}\al_{k}^{2}<\infty$, so that $\calLP'\subset \calLP$. For 
$m\in \bZ_+$ we set $\calLP'_{m}=\calLP'\cap \calLP_m$. By taking constant 
sequences $\bs=\{s\}_{i=0}^{\infty}$ and 
$\bt=\{t\}_{i=0}^{\infty}$ in Corollary~\ref{monot} we obtain 
the following generalization of  Theorems 1.4 and 1.6 in \cite{BS}.

\begin{corollary}\label{coromonot1}
Let $n\in \bN$ and $\vf\in \calLP'$ with $\vf(0)=1$. If $s,t\in \bR$ are 
such that $|s|\le |t|$ and $st\ge 0$ then $\vf(sD)[P]\pc \vf(tD)[P]$ for any 
$P\in \calH_n$.\hfill$\square$ 
\end{corollary}

Let $\calA'$ be the submonoid of $\calA$ 
consisting of all operators that preserve the barycenter of the 
zeros of any nonconstant polynomial. Then by~\eqref{mono2} one has
\begin{equation*}
\begin{split}
\calA' & =\left\{T\in \calA\mid \fm\big(T(P)\big)=\fm(P)\text{ if } P\in \vP,
\,\deg P\ge 1\right\}\\
& = \left\{\vf(D)\mid \vf\in \calLP',\,\vf(0)=1\right\}\subset \calLP'_{0}.
\end{split}
\end{equation*}
Setting $s=0$ and $t=1$ in Corollary~\ref{coromonot1} we deduce that any 
nonconstant monic hyperbolic polynomial is the global 
minimum of its $\calA'$-orbit. In this way we recover Theorem 6 of \cite{B}:

\begin{corollary}\label{coromonot2}
If $n\in \bN$ then $P\pc T[P]$ for all $P\in \calH_n$ and 
$T\in \calA'$.\hfill$\square$
\end{corollary}

Finally, let us note that some of the properties established above may be 
restated by using the following terminology of set-theoretic topology:

\begin{definition}\label{d-set}
An operator $T$ on a poset $(\calX,\le)$ is called {\em isotonic} if 
$T[x]\le T[y]$ whenever $x,y\in \calX$ are such that $x\le y$ while 
$T$ is said to be {\em extensive} (or {\em expanding}) if 
$x\le T[x]$ for any $x\in \calX$. An operator on $(\calX,\le)$ which 
is idempotent, isotonic and extensive with respect to 
$\le$ is called a {\em closure operator} on $\calX$.
\end{definition}

For instance, Corollary~\ref{iso} asserts that 
essentially all differential operators of Laguerre-P\'olya type are isotonic 
on each of the posets $(\calH_n,\pc)$, $n\in \bN$, while 
Corollary~\ref{coromonot1} shows that the monoid $\calA'$ 
consists of differential operators of Laguerre-P\'olya type which are 
extensive with respect to the spectral order. 

\begin{remark}
The proofs of Theorems~\ref{main1} and~\ref{main2} were essentially
based on a detailed analysis of the dynamics of the zeros and critical 
points of strictly hyperbolic polynomials under the action of differential 
operators of Laguerre-P\'olya type. There are many known examples of such 
operators that actually map any hyperbolic polynomial to a {\em strictly} 
hyperbolic polynomial (cf., e.g., \cite{CC1,CC2}). For instance, if $Q$ is a 
hyperbolic polynomial of 
degree $n$ and $b\in \bR$ then $e^{bD}Q(D)[P]$ is strictly hyperbolic 
whenever $P$ is a hyperbolic polynomial of degree at most $n+1$. Moreover, if 
$\vf(x)$ is a transcendental function in the Laguerre-P\'olya class which is 
not of the form $Q(x)e^{bx}$ for some hyperbolic 
polynomial $Q$ and $b\in \bR$ then a theorem of P\'olya asserts that 
$\vf(D)[P]$ is strictly hyperbolic for any hyperbolic polynomial $P$. In 
particular, this holds if $\vf(x)=e^{-a^{2}x}$ with $a\in \bR\setminus \{0\}$. 
\end{remark}

\section{Further results and related topics}\label{s4}

In this section we state several other consequences of Theorems~\ref{main1} 
and~\ref{main2} and discuss some related problems.

\subsection{The distribution of zeros of hyperbolic polynomials}

The results given in \S \ref{s2}--\ref{s3} have interesting applications to
the distribution and the relative geometry of zeros of hyperbolic polynomials 
and their images under the action of differential operators of 
Laguerre-P\'olya type. Recall from \S \ref{s2} that a function 
$\Phi:\bR^n\rightarrow \bR$ is said to be Schur-convex if $\Phi(X)\le \Phi(Y)$ 
whenever $X,Y\in \bR^n$ are such that $X\pr Y$. Given a polynomial $P\in\vP$ 
of degree $n\ge 1$ we denote its zeros by $x_i(P)$, $1\le i\le n$. Then
Theorems~\ref{spec} and Corollary~\ref{iso} yield the following result.

\begin{corollary}\label{con1}
Let $n,m\in \bZ_+$ with $n\ge m+1$. If $\vf\in\calLP_m$ and 
$\Phi:\bR^{n-m}\rightarrow \bR$ is a Schur-convex function then
$$\Phi\big(x_1(\vf(D)[Q]),\ldots,x_{n-m}(\vf(D)[Q])\big)\le 
\Phi\big(x_1(\vf(D)[P]),\ldots,x_{n-m}(\vf(D)[P])\big)$$
for all polynomials $P,Q\in \calH_n$ such that $Q\pc P$. In particular, the 
inequality
$$\sum_{i=1}^{n-m}f\big(x_i(\vf(D)[Q])\big)\le 
\sum_{i=1}^{n-m}f\big(x_i(\vf(D)[P])\big)$$
holds for any convex function $f:\bR\rightarrow \bR$.\hfill$\square$
\end{corollary}

In the same spirit, Theorem~\ref{main2} and 
Corollaries~\ref{coromonot1}--\ref{coromonot2} combined with 
Theorem~\ref{spec} lead to the following inequalities.

\begin{corollary}\label{con2}
Let $n\in \bN$ and $\vf\in \calLP_{0}'$. For any pair 
$(s,t)\in \bR^2$ satisfying $|s|\le |t|$ and $st\ge 0$ and for any
Schur-convex function $\Phi:\bR^{n}\rightarrow \bR$ one has
$$\Phi\big(x_1(\vf(sD)[P]),\ldots,x_{n}(\vf(sD)[P])\big)\le 
\Phi\big(x_1(\vf(tD)[P]),\ldots,x_{n}(\vf(tD)[P])\big)$$
whenever $P\in \calH_n$. In particular, the inequalities
\begin{equation*}
\begin{split}
\sum_{i=1}^{n}f\big(x_i(\vf(sD)[P])\big) & \le 
\sum_{i=1}^{n}f\big(x_i(\vf(tD)[P])\big)\\
\sum_{i=1}^{n}f\big(x_i(P)\big) &\le 
\sum_{i=1}^{n}f\big(x_i(\vf(tD)[P])\big)
\end{split}
\end{equation*}
hold for any convex function $f:\bR\rightarrow \bR$.\hfill$\square$
\end{corollary}

Let $\calLP''$ denote the class of entire functions of the form
$$\vf(x)=cx^{m}e^{bx}\prod_{k=1}^{\infty}(1-\al_k x),$$
where $c\in\bR\setminus \{0\}$, $m\in \bZ_+$, $b\le 0$, $\al_k\ge 0$ and 
$\sum_{k=1}^{\infty}\al_k<\infty$, so that $\calLP''\subset \calLP$. It is 
well-known that $\calLP''$ consists precisely of those functions which are 
locally uniform limits in $\bC$ of sequences of hyperbolic polynomials 
having only positive zeros (cf.~\cite[Ch.~8]{L}). According to the terminology 
introduced by P\'olya and Schur, a real entire function $\psi$ is called a 
function of {\em type I} in the Laguerre-P\'olya class if either 
$\psi(x)\in\calLP''$ or $\psi(-x)\in\calLP''$.
For $m\in \bZ_+$ we set $\calLP_{m}''=\calLP''\cap\calLP_m$. Let 
$P\in \calH_n$ with $n\ge 1$ be such that $x_i(P)>0$ for $1\le i\le n$. Using 
Lemma~\ref{l1} and polynomial approximations as in~\eqref{approx} 
and~\eqref{lim1} one can show that if $\vf\in\calLP_{m}''$ and $n\ge m+1$ 
then $x_i(\vf(D)[P])>0$ for $1\le i\le n-m$. These observations allow us to 
derive new inequalities involving differential operators associated with 
functions of type I in the Laguerre-P\'olya class. The first two inequalities 
listed in Corollary~\ref{con3} below correspond to the following special 
choices of convex functions in Corollary~\ref{con1}: 
minus the Shannon entropy $-H(x)=x\log x$ and minus the Renyi entropies 
$\log(\sum_{i=1}^{n}x_{i}^k)$ for $k\ge 1$, respectively. These are in fact
easy consequences of the third inequality, which is actually the most 
general inequality of this type.

\begin{corollary}\label{con3}
Let $n,m\in \bZ_+$ with $n\ge m+1$. For any $\vf\in\calLP_{m}''$ one has
\begin{equation*}
\begin{split}
\sum_{i=1}^{n-m}x_i(\vf(D)[Q])\log x_i(\vf(D)[Q]) & \le 
\sum_{i=1}^{n-m}x_i(\vf(D)[P])\log x_i(\vf(D)[P]),\\
\sum_{i=1}^{n-m}[x_i(\vf(D)[Q])]^k & \le \sum_{i=1}^{n-m}[x_i(\vf(D)[P])]^k,
\quad k\in [1,\infty),\\
r(r-1)\sum_{i=1}^{n-m}[x_i(\vf(D)[Q])]^r & \le 
r(r-1)\sum_{i=1}^{n-m}[x_i(\vf(D)[P])]^r,\quad r\in \bR,
\end{split}
\end{equation*}
for all polynomials $P,Q\in \calH_n$ with positive zeros that satisfy
$Q\pc P$.\hfill$\square$
\end{corollary}

\subsection{Multiplier sequences, spectral order and isotonic operators}

It is na\-tural to ask whether the spectral order is preserved by linear 
operators other than those of Laguerre-P\'olya 
type (cf.~Problem~\ref{pb3} below). Clearly, any such operator should 
necessarily map 
hyperbolic polynomials to hyperbolic polynomials of the same degree. An 
important class of operators that one may consider in this context is 
the class of diagonal operators (in the basis of standard monomials) that 
preserve hyperbolicity. This is the class of multiplier sequences of the 
first kind, which was completely characterized by P\'olya and Schur 
in~\cite{PS}. 

\begin{definition}\label{ms}
Let $\Ga=\{\ga_k\}_{k=0}^{\infty}$ be an arbitrary sequence of real numbers
and let $T_{_{\Ga}}\in \E\vP$ be given by 
$T_{_{\Ga}}[x^n]=\ga_{n}x^n$, $n\in \bZ_+$. Then $\Ga$ is called 
a {\em multiplier sequence of the first kind} if 
$T_{_{\Ga}}$ preserves the class of hyperbolic polynomials.
\end{definition}

For convenience, we denote by $\calPS_{I}$ the set of all multiplier sequences 
of the first kind and we let $\vP_n$ be the 
$(n+1)$-dimensional subspace of $\vP$ consisting of all complex polynomials 
of degree at most $n$, so that $\calH_n\subset \vP_n$. If 
$\Ga=\{\ga_k\}_{k=0}^{\infty}\in \calPS_{I}$ and $\ga_n\neq 0$ for 
some $n\in \bN$ we define the 
{\em $n$-th normalized truncation} of $\Ga$ to be the finite sequence 
$\Ga_n=\left\{\dfrac{\ga_0}{\ga_n},\ldots,\dfrac{\ga_{n-1}}{\ga_n},1\right\}$. 
Obviously, $\Ga_n$ induces a well-defined 
linear operator $T_{_{\Ga_n}}\in \E\vP_n$ that satisfies 
$T_{_{\Ga_n}}(\calH_n)\subseteq \calH_n$.

\begin{problem}\label{pb1}
Let $\Ga=\{\ga_k\}_{k=0}^{\infty}\in \calPS_{I}$ be such that 
$\ga_n\neq 0$, $n\in \bN$. Is it true that for any $n\in \bN$ the  
operator $T_{_{\Ga_n}}\in \E\vP_n$ preserves the partial ordering $\pc$ on 
$\calH_n$, where $\Ga_n$ is the $n$-th normalized truncation of $\Ga$?
\end{problem}

The condition $\ga_n\neq 0$, $n\in \bN$, 
imposed in Problem~\ref{pb1} is far from being as restrictive as it may 
first appear and is actually quite natural in view of well-known 
properties of multiplier sequences of the first kind (see, e.g., \cite{L}). 
Indeed, if $\Ga=\{\ga_k\}_{k=0}^{\infty}\in \calPS_{I}$ then 
$\{\ga_{i+k}\}_{k=0}^{\infty}\in \calPS_{I}$ for any $i\in \bN$. Moreover, 
if $\ga_0\neq 0$ and $\ga_i=0$ for some $i\in \bN$ then $\ga_j=0$ for all 
$j\ge i$. It follows that either $\Ga$ contains only zero terms except for 
a finite number of consecutive nonzero elements or there exists $i\in \bZ_+$ 
such that $\ga_k=0$ for $k\le i-1$ and $\ga_k\neq 0$ if $k\ge i$. 

As an example, consider the sequence $\Ga=\{k\}_{k=0}^{\infty}$ consisting of 
the Maclaurin coefficients of $xe^x$. Clearly, 
$T_{_{\Ga}}[P(x)]=xP'(x)$ for any $P\in \vP$ hence 
$T_{_{\Ga_n}}(\calH_n)\subseteq \calH_n$, $n\in \bN$. Note that in this 
case Lemma~\ref{deriv} and Theorem~\ref{spec} imply that 
$T_{_{\Ga_n}}$ preserves indeed all the poset 
structures $(\calH_n,\pc)$, $n\in \bN$ . Similar considerations show that the 
answer to Problem~\ref{pb1} is affirmative for multiplier sequences of the 
following type.

\begin{proposition}\label{lag-ms}
Let $m\in \bN$, $p\in \bZ_+$ and consider the sequence 
$\Ga=\{H(k+p)\}_{k=0}^{\infty}$, where
$H(x)=\prod_{i=0}^{m-1}(x-i)$. Then $\Ga\in \calPS_{I}$ and for 
any $n\ge \max(1,m-p)$ the operator 
$T_{_{\Ga_n}}$ preserves the partial ordering $\pc$ on $\calH_n$.
\end{proposition}

\begin{proof}
If $n\in \bN$ and $P(x)=\sum_{k=0}^{n}x^k\in \vP_n$ then
$$T_{_{\Ga}}[P(x)]=\sum_{k=0}^{n}H(k+p)a_k x^k
=x^{m-p}\left[x^{p}P(x)\right]^{(m)}$$
and so by Rolle's theorem $\Ga$ is a multiplier 
sequence of the first kind. The same arguments further show that  
$T_{_{\Ga_n}}(\calH_n)\subseteq \calH_n$ for all $n\ge \max(1,m-p)$ 
since $H(n+p)\neq 0$ for such $n$. Using Lemma~\ref{deriv}
and Theorem~\ref{spec} (i) one can easily check that 
$x^{m-p}\left[x^{p}Q(x)\right]^{(m)}\pc x^{m-p}\left[x^{p}P(x)\right]^{(m)}$ 
whenever $n\ge \max(1,m-p)$ and $P,Q\in \calH_n$ are such that $Q\pc P$.
\end{proof}

A somewhat different version of Problem~\ref{pb1} is as follows.

\begin{problem}\label{pb2}
Fix $n\in \bN$ and consider a finite sequence 
$\La=\{\la_k\}_{k=0}^{n}$ with associated operator 
$T_{_{\La_n}}\in \E\vP$ given by $T_{_{\La_n}}[x^k]=\la_k x^k$,
$0\le k\le n$, $T_{_{\La_n}}[x^k]=0$, $k>n$. If $\la_n=1$ and 
$T_{_{\La_n}}(\calH_n)\subseteq \calH_n$ 
is it true that $T_{_{\La_n}}$ preserves the spectral order on $\calH_n$?
\end{problem}

The answer to Problem~\ref{pb2} is trivially affirmative if $n=1$ and 
elementary computations show that this holds for $n=2$ as well. Indeed, if 
$\La=\{\la_0,\la_1,1\}$ is a sequence that satisfies the above hypotheses 
then $\la_0\ge 0$ since $T_{_{\La_n}}[x^2-1]\in \calH_2$. Given two 
polynomials $P(x)=x^2+ax+b\in \calH_2$ and $Q(x)=x^2+cx+d\in \calH_2$ with 
$Q\pc P$ one has $a=c$, $a^2\ge 4\max(b,d)$ and 
$\sqrt{a^2-4d}\le \sqrt{a^2-4b}$. From $\la_0\ge 0$ we get
$\sqrt{\la_{1}^{2}a^2-4\la_{0}d}\le \sqrt{\la_{1}^{2}a^2-4\la_{0}b}$, which 
shows that $T_{_{\La_n}}[Q]\pc T_{_{\La_n}}[P]$.

Problem~\ref{pb2} may actually be viewed as a special case of a yet more 
general problem. Fix $n\in \bN$ and recall the monoid $\calA_n$ defined 
in~\eqref{mono1}. Let 
$\calA_{n}^{\pc}$ denote the submonoid of $\calA_n$ consisting of all 
operators that preserve the poset structure $(\calH_n,\pc)$, that is,
$$\calA_{n}^{\pc}
=\{T\in \calA_n\mid T[Q]\pc T[P]\text{ if }P,Q\in \calH_n,\,Q\pc P\}.$$
Recall also the submanifold $\calH_{n}^{0}$ of $\calH_n$ from~\eqref{h0} 
and consider the submonoid $\calA_{n}^{0}$ of $\calA_n$ given by
$$\calA_{n}^{0}=\{T\in \calA_n\mid T\big(\calH_{n}^{0}\big)\subseteq 
\calH_{n}^{0}\big\}.$$

\begin{problem}\label{pb3}
Describe all operators in $\calA_{n}^{\pc}$. Is it true that 
$\calA_{n}^{\pc}=\calA_{n}^{0}$ for all $n\in \bN$?
\end{problem}

\begin{conjecture}\label{c-true}
Problems~\ref{pb1}--\ref{pb3} have all affirmative answers.
\end{conjecture}

\begin{remark}
The linear 
transformations on $\bR^n$ that preserve the majorization relation $\pr$ 
between $n$-vectors of real numbers were characterized in~\cite{An2,DV}.
\end{remark}

Note that Problem~\ref{pb3} implicitly addresses and further motivates both 
the question of describing all operators in the monoid $\calA_n$ itself 
(cf.~\cite[Problem 2 (iii)]{B}) and its version with no restriction on 
the degrees that may be formulated as follows.

\begin{problem}\label{pb4}
Characterize all operators in the 
monoid $\tilde\calA:=\{T\in \E\vP\mid T(\calH)\subseteq \calH\}$, where 
$\calH=\bigcup_{n=0}^{\infty}\calH_n$.
\end{problem}

Problem~\ref{pb4} is actually a long-standing open problem of 
fundamental interest in the theory of distribution of zeros of polynomials 
and transcendental entire functions (see \cite[Problem 1.3]{CC1}). 
Significant progress towards a complete
solution to Problem~\ref{pb4} was recently made in \cite{BBS}.

The above results and those of~\cite{B,BP,BS} show that even 
a partial knowledge of the operators in $\calA_n$ leads to new 
interesting information on 
the relative geometry of the zeros of a hyperbolic polynomial and the zeros of 
its images under such operators. Several related questions arise naturally in 
this context. For instance, Problem 2 (ii) in~\cite{B} asks whether it is 
possible to describe the spectral order by means of the action of linear 
(differential) operators on the partially ordered mani\-fold $(\calH_n,\pc)$. 
This would provide a new characterization of classical majorization which in 
a way would be dual to the usual characterization in terms of doubly 
stochastic matrices given in Theorem~\ref{spec}.

It would also be 
interesting to know whether there are any ``infinite-dimensional'' analogs 
of Theorems~\ref{main1} and~\ref{main2}. Indeed, it is well known that the 
class $\calLP$ is closed under differentiation \cite{L}. A more 
general closure property was established in~\cite{CC2}, where various 
types of infinite order differential operators acting on 
$\calLP$ were studied in detail. In particular, Lemmas 3.1 and 3.2 in 
{\em loc.~cit.~}show that the subset of $\calLP$ consisting of 
entire functions of genus 0 or 1 
is stable under the action of differential operators 
of Laguerre-P\'olya type. Moreover, there are several known extensions of 
classical majorization to infinite sequences of real numbers 
\cite[p.~16]{MO}. One may therefore ask if these extensions or some 
appropriate modifications could lead to 
generalizations of the above results to differential operators acting on 
transcendental entire functions in the class $\calLP$.

\end{document}